\documentclass[12pt,leqno]{article}
\usepackage{amsmath,amsfonts}
\usepackage[dvips]{graphics,color}
\usepackage{latexsym}
\usepackage{amssymb}
\usepackage{graphicx}
\date{ }
\begin{document}
\title{ Embedded constant $m^{\text{th}}$ mean curvature hypersurfaces on spheres
\thanks{The first author was supported
by grant No. 10971110 of NSFC and by the Doctoral Program Foundation of the
Ministry of Education of China (Grant No. 20104407120002).}}
\author{Guoxin Wei and Guohua Wen } \maketitle

\begin{abstract}
\noindent In this paper, we study $n$-dimensional hypersurfaces with
constant $m^{\text{th}}$ mean curvature $H_m$ in a unit sphere
$S^{n+1}(1)$ and prove that if the
$m^{\text{th}}$ mean curvature $H_m$ takes value between
 $\dfrac{1}{(\tan \frac{\pi}{k})^m}$ and $\frac{k^2-2}{n}(\frac{k^2+m-2}{n-m})^{\frac{m-2}{2}}$ for $1\leq m\leq n-1$
 and any integer $k\geq 2$, then there exists at least one
$n$-dimensional compact nontrivial embedded hypersurface with
constant $H_m>0$ in $S^{n+1}(1)$. When $m=1$, our results reduce to the results of Perdomo \cite{[P]}; when $m=2$ and $m=4$, our results reduce to the results of Cheng-Li-Wei \cite{[WCL]}.

\end{abstract}



\section *{1. Introduction}

It is well known that  Alexandrov \cite{[A]} and  Montiel-Ros
\cite{[MR]} proved that the standard round spheres are the only
possible oriented compact embedded hypersurfaces with constant
$m^{\text{th}}$ mean curvature $H_m$ in a Euclidean space
$\mathbb{R}^{n+1}$ , for $m\geq 1$. For hypersurfaces in a unit sphere $S^{n+1}(1)$,
standard round spheres and Clifford hypersurfaces
$S^{l}(a)\times S^{n-l}(b)$, $1\leq l\leq n-1$ are compact embedded
hypersurfaces in $S^{n+1}(1)$. Hence the following problem is interesting (also see \cite{[BL]}, \cite{[Le]}, \cite{[WCL]}):

\bigskip \noindent
{\bf Problem:} Do there exist  compact  embedded hypersurfaces with constant
$m^{\text{th}}$ mean curvature $H_m$ in  $S^{n+1}(1)$ other than the
standard round spheres and Clifford
hypersurfaces?\\

\noindent When $m=1$, namely, when the mean curvature is constant, Ripoll
\cite{[R]} has proved the existence of compact embedded
hypersurfaces of $S^3(1)$ with constant mean curvature ($H\neq0,
\pm\frac{\sqrt{3}}{3}$) other than the standard round spheres and the
Clifford hypersurfaces.  For general $n$, Perdomo \cite{[P]} has proved

\par\bigskip \noindent{\bf Theorem 1.1} (Main Theorem of \cite{[P]}). {\it For any $n\geq 2$ and any
integer $k\geq 2$, if mean curvature $H$ takes
value between $\dfrac{1}{(\tan \frac{\pi}{k})}$ and
$\dfrac{(k^2-2)\sqrt{n-1}}{n\sqrt{k^2-1}}$, then there exists an $n$-dimensional
compact nontrivial embedded hypersurface with constant mean curvature $H>0$ in
$S^{n+1}(1)$.}\\

\noindent For $m=2$, that is, when the scalar curvature is constant, Cheng, Li and Wei \cite{[WCL]} has proved
\par\bigskip \noindent{\bf Theorem 1.2} (\cite{[WCL]}). {\it For any $n\geq3$ and any
integer $k\geq 2$, if $H_2=\frac{R-n(n-1)}{n(n-1)}$ takes value
between $\dfrac{1}{(\tan \frac{\pi}{k})^2}$ and $\dfrac{k^2-2}{n}$,
then there exists an $n$-dimensional compact nontrivial embedded
hypersurface $M$ with constant 2-th mean curvature $H_2>0$ $($i.e.
scalar curvature $R>n(n-1)$$)$ in  $S^{n+1}(1)$, where $R$ is the
scalar curvature of $M$.}\\

\noindent For $m=4$, Cheng, Li and Wei \cite{[WCL]} has proved
\par\bigskip \noindent{\bf Theorem 1.3} (\cite{[WCL]}). {\it For any $n\geq 5$ and any
integer $k\geq 3$, if $4^{\text{th}}$ mean curvature $H_4$ takes
value between $\dfrac{1}{(\tan \frac{\pi}{k})^4}$ and
$\dfrac{k^4-4}{n(n-4)}$, then there exists an $n$-dimensional
compact nontrivial embedded hypersurface with constant $H_4>0$ in
$S^{n+1}(1)$.}\\

\noindent For general $1\leq m\leq n-1$, we will prove that there exist many compact
nontrivial embedded hypersurfaces with constant  $m^{\text{th}}$
mean curvature $H_m>0$  in $S^{n+1}(1)$, for  $1\leq m\leq n-1$. In fact, we prove

\par\bigskip \noindent{\bf Theorem 1.4}. {\it For $1\leq m\leq n-1$ and any
integer $k\geq 2$, if $m^{\text{th}}$ mean curvature $H_m$ takes
value between $\dfrac{1}{(\tan \frac{\pi}{k})^m}$ and $\dfrac{k^2-2}{n}\biggl(\dfrac{k^2+m-2}{n-m}\biggl)^{\frac{m-2}{2}}$, then there exists at least one $n$-dimensional
compact nontrivial embedded hypersurface with constant $H_m>0$ in
$S^{n+1}(1)$.}

\par\bigskip \noindent{\bf Remark 1.1}. {\it For $m=1$, Theorem 1.4 reduces to the conclusion of Perdomo \cite{[P]}.
For $m=2$, $m=4$, Theorem 1.4 reduces to the results of Cheng-Li-Wei \cite{[WCL]}}.

\section*{2. Proof of Theorem}

Using the same notations as those of \cite{[WCL]}, we can have
 $$(g^{'})^2=q(g),\ \ \text{where}\ \
 q(v)=C-v^2(v^{-n}+H_m)^{\frac{2}{m}}-v^2,\eqno{(2.1)}$$

\begin{equation*}
\aligned &\ \ \ q^{'}(v)\\
&=2v\left\{-(v^{-n}+H_m)^{\frac{2}{m}}+\frac{n}{m}v^{-n}(v^{-n}+H_m)^{\frac{2-m}{m}}-1\right\}\\
&=-2v\left\{(v^{-n}+H_m)^{\frac{2-m}{m}}\left[\frac{m-n}{m}v^{-n}+H_m\right]+1\right\},
\endaligned
\eqno{(2.2)}
\end{equation*}
and
\begin{equation*}
\aligned &\ \ \ q^{''}(v)\\
&=-\frac{2(v^{-n}+H_m)^{\frac{2-2m}{m}}}{m^2}\left\{(2n^2-3nm+m^2)v^{-2n}+m(n^2-3n+2m)H_mv^{-n}+m^2H_m^2\right\}\\
&\ \ \ -2\\
&<-2.
\endaligned
\eqno{(2.3)}
\end{equation*}
From (2.3), one obtains that $q^{'}(v)$ is a decreasing function of $v$ in $[0, +\infty)$. From (2.2), one has
 $q^{'}(v)>0$ if $v\rightarrow 0$; $q^{'}(v)<0$ if $v\rightarrow \infty$. Hence there exists $0<v_0<\infty$ such that
 $q^{'}(v_0)=0$. Moreover, the function $q(v)$ is a monotone increasing function of $v$
in $(0, v_0]$ and decreasing function of $v$ in $[v_0, +\infty)$. Hence, for some value of $C$, the function $q$
has two positive roots $t_1$ and $t_2$, such that $t_1\leq t_2$, $q(t_1)=q(t_2)=0$ and $q(t)>0$ if $t\in (t_1,t_2)$.

From the results of section 4 of \cite{[WCL]}, we have
$$P(H_m,n,C)=2\int_0^{\frac{T}{2}}\frac{r(s)\lambda(s)}{1-r^2(s)}ds.\eqno{(2.4)}$$
Since
$r(s)=\frac{g(s)}{\sqrt{C}}$ and
$\lambda(s)=(g(s)^{-n}+H_m)^{\frac{1}{m}}$, then it follows from
(2.4) that
$$P(H_m,n,C)=2\int_0^{\frac{T}{2}}\frac{\sqrt{C}g(s)(g(s)^{-n}+H_m)^{\frac{1}{m}}}{C-g^2(s)}ds.\eqno{(2.5)}$$
Doing the substitution $t=g(s)$, putting $c_0=C-q(v_0),\ \ -2a=q^{''}(v_0)$ and applying Lemma 5.1 of \cite{[WCL]}, one concludes from $g(0)=t_1$ and $g(\frac{T}{2})=t_2$ that

\begin{equation*}
\aligned &\lim_{C\rightarrow c_0^{+}}
P(H_m,n,C)\\
 &=\frac{2\pi\sqrt{c_0}}{\sqrt{a}\sqrt{c_0-v_0^2}}\\
 &=\dfrac{2m\pi\biggl((v_0^{-n}+H_m)^{\frac{2m-2}{m}}+(v_0^{-n}+H_m)^{\frac{2m-4}{m}}\biggl)^{\frac{1}{2}}}
 {\biggl((2n^2-3mn+m^2)v_0^{-2n}+m(n^2-3n+2m)H_mv_0^{-n}+m^2H_m^2+m^2(v_0^{-n}+H_m)^{\frac{2m-2}{m}}\biggl)^{\frac{1}{2}}}.
\endaligned
\eqno{(2.6)}
\end{equation*}
Writing
$$F_0=(v_0^{-n}+H_m)^{\frac{1}{m}}, \eqno{(2.7)}$$
 then $F_0>0$ and $q^{'}(v_0)=0$ can be written as
$$F_0^m+\frac{m}{m-n}F_0^{m-2}+\frac{n}{m-n}H_m=0.\eqno{(2.8)}$$
We next compute the the numerator and denominator of $P(H_m,n,C)$ by using of (2.8). By a direct calculation, we know
\begin{equation*}
\aligned
&2m\pi\biggl((v_0^{-n}+H_m)^{\frac{2m-2}{m}}+(v_0^{-n}+H_m)^{\frac{2m-4}{m}}\biggl)^{\frac{1}{2}}\\
&=2m\pi\Biggl(F_0^{m-2}\biggl(\dfrac{m}{n-m}F_0^{m-2}+\frac{n}{n-m}H_m\biggl)+F_0^{2m-4}\Biggl)^{\frac{1}{2}}\\
&=2m\pi\biggl(\frac{n}{n-m}\biggl)^{\frac{1}{2}}F_0^{\frac{m-2}{2}}\biggl(F_0^{m-2}+H_m\biggl)^{\frac{1}{2}},\\
\endaligned
\eqno{(2.9)}
\end{equation*}

\begin{equation*}
\aligned
&\ (2n^2-3mn+m^2)v_0^{-2n}+m(n^2-3n+2m)H_mv_0^{-n}+m^2H_m^2\\
&\ \ \ +m^2(v_0^{-n}+H_m)^{\frac{2m-2}{m}}\\
&=(n-m)(2n-m)\biggl(v_0^{-n}+H_m\biggl)^2+n(-4n+mn+3m)(v_0^{-n}+H_m)H_m\\
&\ \ \ +n^2(2-m)H_m^2+m^2(v_0^{-n}+H_m)^{\frac{2m-2}{m}}\\
&=(n-m)(2n-m)\biggl(\frac{m}{n-m}F_0^{m-2}+\frac{n}{n-m}H_m\biggl)^2\\
&\ \ \ +n(-4n+mn+3m)\biggl(\frac{m}{n-m}F_0^{m-2}+\frac{n}{n-m}H_m\biggl)H_m
+n^2(2-m)H_m^2\\
&\ \ \ +m^2F_0^{m-2}\biggl(\frac{m}{n-m}F_0^{m-2}+\frac{n}{n-m}H_m\biggl)^{\frac{2m-2}{m}}\\
&=\frac{m^2n}{n-m}\biggl(F_0^{m-2}+H_m\biggl)\biggl(2F_0^{m-2}+nH_m\biggl).\\
\endaligned
\eqno{(2.10)}
\end{equation*}
According (2.6), (2.8), (2.9) and (2.10), we obtain
\begin{equation*}
\aligned &\lim_{C\rightarrow c_0^{+}}
P(H_m,n,C)\\
 &=\dfrac{2\pi\sqrt{c_0}}{\sqrt{a}\sqrt{c_0-v_0^2}}\\
 &=\dfrac{2m\pi\biggl((v_0^{-n}+H_m)^{\frac{2m-2}{m}}+(v_0^{-n}+H_m)^{\frac{2m-4}{m}}\biggl)^{\frac{1}{2}}}
 {\biggl((2n^2-3mn+m^2)v_0^{-2n}+m(n^2-3n+2m)H_mv_0^{-n}+m^2H_m^2+m^2(v_0^{-n}+H_m)^{\frac{2m-2}{m}}\biggl)^{\frac{1}{2}}}\\
 &=2\pi\biggl( \dfrac{F_0^{m-2}}{(n-m)F_0^m-(m-2)F_0^{m-2}}\biggl)^{\frac{1}{2}}\\
 &=2\pi\biggl( \dfrac{1}{(n-m)F_0^2-(m-2)}\biggl)^{\frac{1}{2}},
\endaligned
\eqno{(2.11)}
\end{equation*}
and
$$F_0^2>\frac{m-2}{n-m}.\eqno{(2.12)}$$

On the other hand, we know that
$$\lim_{C\rightarrow
\infty}P(H_m,n,C)=2\arctan\frac{1}{H_m^{\frac{1}{m}}}.\eqno{(2.13)}$$
Therefore, for any fixed $H_m>0$, the function $P(H_m,n,C)$ takes
all the values between
$$A(H_m)=2\arctan\frac{1}{H_m^{\frac{1}{m}}},\ \
B(H_m)=2\pi\biggl( \frac{1}{(n-m)F_0^2-(m-2)}\biggl)^{\frac{1}{2}}. \eqno{(2.14)}$$
By a direct calculation, we obtain $A(H_m)$ is a decreasing function of $H_m$. From (2.8) and (2.12), it is not hard to prove $F_0$ is an increasing function of $H_m$, then $\lim\limits_{C\rightarrow c_0^{+}}P(H_m,n,C)$ is a decreasing function of $H_m$. Moreover, if $2\pi\bigl( \dfrac{1}{(n-m)F_0^2-(m-2)}\bigl)^{\frac{1}{2}}=\frac{2\pi}{k}$, then
$$F_0^2=\frac{k^2+m-2}{n-m},\eqno{(2.15)}$$
from (2.8), we obtain
$$H_m=\frac{k^2-2}{n}\biggl(\frac{k^2+m-2}{n-m}\biggl)^{\frac{m-2}{2}},\eqno{(2.16)}$$
thus,
$$A\biggl(\frac{1}{(\tan
\frac{\pi}{k})^m}\biggl)=B\Biggl(  \frac{k^2-2}{n}\biggl(\frac{k^2+m-2}{n-m}\biggl)^{\frac{m-2}{2}}\Biggl)=\frac{2\pi}{k},\eqno{(2.17)}$$
where $k\geq2$ is any integer, then we deduce that the number
$\frac{2\pi}{k}$ lies between $A(H_m)$ and $B(H_m)$. Hence, by the continuity of $P(H_m,n,C)$, there exists
 some constant $C_1$ such that $P(H_m,n,C_1)=\frac{2\pi}{k}$.
If the period is $\frac{2\pi}{k}$, then there exists a compact
embedded hypersurface with constnat $H_m$ which is not isometric to
a round sphere or a Clifford hypersurface. We complete the proof of
Theorem 1.4.$$\eqno{\Box}$$

\begin{flushleft}
\medskip\noindent
\begin{tabbing}
XXXXXXXXXXXXXXXXXXXXXXXXXX*\=\kill Guoxin Wei\> Guohua Wen\\
 School of Mathematical Sciences\> School of Mathematical Sciences \\
 South China Normal University\>South China Normal University\\
  510631, Guangzhou\>510631, Guangzhou\\
 China\>China\\
 E-mail: weigx@scnu.edu.cn\\
 \ \ \ \ \ \ \ \ \ \ \ weigx03@mails.tsinghua.edu.cn \\
\end{tabbing}
\end{flushleft}

\end {document}